\documentclass[
leqno,          
10pt,           
a4paper         
]{amsart}

 \usepackage[colorlinks,citecolor=blue,urlcolor=blue]{hyperref}          
\usepackage{hyperref}
\usepackage{courier}                                                    
\usepackage{amssymb, amsmath, amsfonts, amsthm}                         
\usepackage{mathtools, cite, enumerate, rotating, framed, lipsum}       
\usepackage{fancybox, bbm}                                              
\usepackage[dvipsnames]{xcolor}                                         
\usepackage[polish, english]{babel}                                     
\usepackage[T1]{fontenc}                                                
\usepackage[cp1250]{inputenc}                                           
\DeclareMathAlphabet{\mathpzc}{OT1}{pzc}{m}{it}

\numberwithin{equation}{section}                        

\theoremstyle{definition}

\newcommand{\thmcount}{equation}                 
\newcounter{specialcounter}

\newtheorem{Thm}[\thmcount]{Theorem}
\newtheorem{Sthm}[specialcounter]{Theorem}
\newtheorem{Cor}[\thmcount]{Corollary}
\newtheorem{Lem}[\thmcount]{Lemma}
\newtheorem{Prop}[\thmcount]{Proposition}
\newtheorem{Rem}[\thmcount]{Remark}
\newtheorem{Defn}[\thmcount]{Definition}
\newtheorem{Ex}[\thmcount]{Example}
\newtheorem{Asu}[\thmcount]{Assumption}
\newtheorem{Sol}[\thmcount]{Solution}
\newtheorem{Fact}[\thmcount]{Fact}
\newtheorem*{Thmx}{Theorem}
\newtheorem*{Corx}{Corollary}
\newtheorem*{Lemx}{Lemma}
\newtheorem*{Propx}{Proposition}
\newtheorem*{Remx}{Remark}
\newtheorem*{Defnx}{Definition}
\newtheorem*{Exx}{Example}
\newtheorem*{Asux}{Assumption}
\newtheorem*{Solx}{Solution}
\newtheorem*{Factx}{Fact}

\newcommand \eq[1]{\begin{equation} #1 \end{equation}}
\newcommand \eqx[1]{\begin{equation*}  #1 \end{equation*}}
\newcommand \al[1]{\begin{align} #1 \end{align}}
\newcommand \alx[1]{\begin{align*}  #1 \end{align*}}
\renewcommand \sp[1]{\begin{equation} \begin{split} #1 \end{split} \end{equation}}
\newcommand \spx[1]{\begin{equation*} \begin{split} #1 \end{split} \end{equation*}}

\newcommand{\thm}[2]{\begin{Thm} \label{#1} #2 \end{Thm}}

\newcommand{\lem}[2]{\begin{Lem} \label{#1} #2 \end{Lem}}

\newcommand{\cor}[2]{\begin{Cor} \label{#1} #2 \end{Cor}}

\newcommand{\prop}[2]{\begin{Prop} \label{#1} #2 \end{Prop}}

\newcommand{\defn}[2]{\begin{Defn} \label{#1} #2 \end{Defn}}

\newcommand{\pr}[1]{\begin{proof} #1 \end{proof}}

\newcounter{comcount}

\renewcommand{\hline}{\vbox{\hrule width\textwidth height 1pt}\smallskip}

\renewcommand{\a}{\alpha}       \newcommand{\be}{\beta}         \newcommand{\e}{\varepsilon}
                 \newcommand{\de}{\delta}
        \newcommand{\la}{\lambda}

 \newcommand{\Bb}{\mathbf{B}} 
\newcommand{\CC}{\mathbb{C}}  
  \newcommand{\dd}{\mathcal{D}}

  \newcommand{\hh}{\mathcal{H}}
  \newcommand{\ii}{\mathcal{I}}
  
  \newcommand{\kk}{\mathcal{K}}
  \renewcommand{\ll}{\mathcal{L}}
  
\newcommand{\NN}{\mathbb{N}}

  \newcommand{\qq}{\mathcal{Q}}
\newcommand{\RR}{\mathbb{R}}

 \newcommand{\Uu}{\mathbf{U}} 
  \newcommand{\vv}{\mathcal{V}}

\newcommand{\ZZ}{\mathbb{Z}}  

\newcommand{\supp}{\mathrm{supp}}
\newcommand{\8}{\infty}
\renewcommand{\d}{\partial}
\renewcommand{\le}{\left(}
\renewcommand{\)}{\right)}

\newcommand{\Rd}{{\RR^d}}

\newcommand{\dt}{\frac{\d}{\d t}}

\renewcommand{\rm}[1]{\mathrm{#1}}
\newcommand{\wt}[1]{\widetilde{#1}}
\newcommand{\limn}{\lim_{n\to \8}}
\newcommand{\limk}{\lim_{k\to \8}}

\newcommand{\abs}[1]{\left| #1 \right|}
\newcommand{\set}[1]{\left\{ #1 \right\}}
\newcommand{\norm}[1]{\left\| #1 \right\|}
\newcommand{\expr}[1]{\left( #1 \right)}
\newcommand{\eee}[1]{\left( #1 \right)}

\newcommand{\lab}[1]{\label{#1}}
\newcommand{\sprod}[1]{\langle #1 \rangle}

\newcommand{\st}[1]{\sup_{t\leq#1^2}}
\newcommand{\Pt}[1]{P_{#1}}

\newcommand{\Kt}[1]{K_{#1}}
\newcommand{\opL}{\mathbf{L}} 
\newcommand{\opB}{\mathbf{B}} 

\newcommand{\innprod}[1]{\langle #1 \rangle_\mu}
\newcommand{\vt}{\vartheta}

\newcommand{\ap}{\frac{\a-1}{2}}
\newcommand{\I}[1]{I_{#1}}
\newcommand{\iI}[1]{I_{#1}^{*}}
\newcommand{\iiI}[1]{I_{#1}^{**}}
\newcommand{\iiiI}[1]{I_{#1}^{***}}

\newcommand{\p}[1]{\phi_{#1}}

\newcommand{\Ls}[1]{L^1(#1,\,\mu)}

\newcommand{\ha}[1]{{h}^1_{#1}(\opB)}
\newcommand{\haL}{\mathcal{H}^1(\opL)}

\newcommand{\hatom}[1]{\mathcal{H}_{at}^1(\ii_{#1},\mu)}

\newcommand{\Dom}{\mathrm{Dom}}
\newcommand{\dI}[1]{\I{#1}^{d}}
\newcommand{\xy}{(x,y)}

 \title[ Hardy spaces for Bessel-Schr\"odinger operators ]{ Hardy spaces for Bessel-Schr\"odinger operators }
 \author[ Edyta Kania ]{ Edyta Kania }
 \address{
 Edyta Kania \newline
 \indent Instytut Matematyczny, Uniwersytet Wroc\l awski \newline
 \indent pl. Grunwaldzki 2/4, 50-384 Wroc\l aw, Poland }
 \email{edyta.kania@math.uni.wroc.pl }

 \author[ Marcin Preisner ]{ Marcin Preisner }
 \address{
 Marcin Preisner \newline
 \indent Instytut Matematyczny, Uniwersytet Wroc\l awski \newline
 \indent pl. Grunwaldzki 2/4, 50-384 Wroc\l aw, Poland }
 \email{marcin.preisner@math.uni.wroc.pl }
\subjclass[2010]{42B30, 42B25,  35J10 (primary), 47D03, 43A85 (secondary)}
\thanks{ The second author was supported by Polish funds for sciences grant DEC-2012/05/B/ST1/00672 from Narodowe Centrum Nauki.}
\begin{document}
\begin{abstract}
Consider the Bessel operator with a potential on  $L^2((0,\8), x^\alpha\, dx)$, namely
$$ \mathbf{L} f(x) = -f''(x) - \frac{\alpha}{x}f'(x) +V(x) f(x).$$
We assume that $\a>0$ and $V\in L^1_{loc}((0,\8),x^\alpha\, dx)$ is a non-negative function. By definition, a~function $f\in L^1((0,\8), x^\alpha \, dx)$ belongs to the Hardy space $\mathcal{H}^1(\mathbf{L})$ if
$$\sup_{t>0} \abs{e^{-t\mathbf{L}}f} \in L^1((0,\8), x^\alpha\, dx).$$
Under certain assumptions on $V$ we characterize the space $\mathcal{H}^1(\mathbf{L})$ in terms of atomic decompositions of local type. In the second part we prove that this characterization can be applied to $\mathbf{L}$ for $\alpha \in(0,1)$ with no additional assumptions on the potential $V$.
\end{abstract}
\maketitle                   





\section{Introduction}

\subsection{Background}

The Schr\"odinger operator on $\Rd$ is given by
    \eqx{
    \wt{\ll} f = -\Delta f + \vv\cdot f,
    }
where $\Delta$ is the Laplace operator and $\vv$ is a function called a potential. If we assume that $\vv\in L^1_{loc}(\Rd)$ and $\vv\geq 0$ then one can find a densely defined, self-adjoint operator $\ll$ on $L^2(\Rd)$, that corresponds to $\wt{\ll}$. It is well known that $\ll$ generates the semigroup of contractions $\kk_t = \exp(-t\ll)$ and $\kk_t$ admits an integral kernel $\kk_t(x,y)$ such that
    \alx{
    &\kk_t f(x) = \int_\Rd \kk_t(x,y) f(y) dy,\\
    0\leq &\kk_t(x,y) \leq (4\pi t)^{-d/2} \exp\eee{-\frac{|x-y|^2}{4t}}.
    }

There have been wide studies on harmonic analysis related to Schr\"odinger operators and, more generally, operators with Gaussian bounds. We refer the reader to \cite{Auscher_Ali}, \cite{ADM}, \cite{Bernicot_Zhao}, \cite{Bui}, \cite{CzajaZienkiewicz_ProcAMS}, \cite{Deng}, \cite{Duong_Yan}, \cite{DP_Argentina}, \cite{DP_Arkiv}, \cite{DP_Potential}, \cite{DZ_Studia2}, \cite{DZ_Annali}, \cite{DZ_Revista2}, \cite{DZ_JFAA}, \cite{DZ_Potential_2014}, \cite{Grafakos_Hardy}, \cite{Hofmann_memoirs}, \cite{Shen},  and references therein. In particular, the Hardy spaces
    \eq{\label{defff}
    H^1(\ll) = \set{f \in L^1(\Rd) \ : \ \norm{f}_{H^1(\ll)} := \norm{\sup_{t>0} \abs{\kk_t f} }_{L^1(\Rd)}<\8}
    }
related to $\ll$ were intensively studied. At this point let us mention that the classical Hardy space $H^1(-\Delta)$ has many equivalent definitions, e.g. in terms of: various maximal functions, singular integrals, square functions, etc. A particulary useful result is the atomic decomposition theorem (see \cite{CoifmanWeiss_BullAMS}, \cite{Latter_Studia}): a~function $f\in H^1(-\Delta)$ can be decomposed as $f(x) = \sum_k \la_k a_k(x)$, where $\sum_k |\la_k| \simeq \norm{f}_{H^1(-\Delta)}$ and $a_k$ are classical atoms, that is, there exist balls $B_k$ such that:
    \eqx{\supp \, a_k \subseteq B_k,\quad \norm{a_k}_\8\leq |B_k|^{-1},\quad \int a(x) dx = 0.}
In other words atoms satisfy some localization, size, and cancellation conditions.

Let us mention that Hofmann et al. \cite{Hofmann_memoirs} have found general results (for $\vv$ satisfying $0\leq \vv\in L^1_{loc}(\Rd)$) saying that $H^1(\ll)$ given above is equal to the Hardy spaces via: square functions, atomic or molecular decompositions. However, atoms used in \cite{Hofmann_memoirs} are given in terms of $\ll$, to be more precise: a function $a\in L^2(\Rd)$ is an atom if there exist a ball $B$ and $b \in \mathrm{Dom}(\ll)$, such that: $a=\ll b$, $\supp \, b \subseteq B$ and $b, \ll b$ satisfy some size condition.

An another approach, started by Dziuba\'nski and Zienkiewicz in the 90's, was to find atomic spaces with simple geometric conditions that characterize $H^1(\ll)$. It appeared that this cannot be done in full generality, and~the properties of atoms depend strictly on the potential $\vv$ and the dimension~$d$. For example, if $d\geq 3$, $\vv\in C_c^\8(\Rd)$, and $\vv\not\equiv 0$, then atoms have modified cancellation condition $\int a(x) \omega(x) \, dx=0$, where $\omega$~is such that $0<C^{-1}\leq \omega(x) \leq C$. For this result and generalizations see \cite{DZ_Annali}, \cite{DZ_JFAA}, \cite{Preisner2}. Other results, see \cite{DZ_Studia2}, lead to Hardy spaces with local atoms. It was first observed by Goldberg \cite{Goldberg_Duke} that if we take supremum for $0<t\leq \tau^2$ in \eqref{defff}, then one obtains atomic space with classical atoms complemented with the atoms of~the form $|B|^{-1}\mathbbm{1}_B(x)$, where the ball $B$ has radius $\tau$. In \cite{DZ_Studia2} the authors assume that for $0\leq \vv \in L^1_{loc}(\Rd)$ there exists a family of cubes $\qq = \set{Q_k \ : \ k\in \NN}$ such that
    \eqx{
    \overline{\cup_k Q_k} = \Rd, \quad |Q_k\cap Q_j|=0 \text{ for }k\neq j, \quad d(Q_k)\simeq d(Q_j) \text{ if } Q_k^{***}\cap Q_j^{***} \neq \emptyset.
    }
Here $d(Q)$ is the diameter of $Q$ and $Q^*$ is a cube that has the same center as $Q$ but with slightly enlarged diameter. The atomic space $H^1_{at}(\qq)$ is built on classical atoms and atoms of the form $|Q_k|^{-1} \mathbbm{1}_{Q_k}(x)$. The main result of \cite{DZ_Studia2} states that under two additional assumptions on $V,\qq,\kk_t$ (see \cite[p.41]{DZ_Studia2}, conditions: $(D)$, $(K)$)  we have that $H^1(\ll) = H^1_{at}(\qq)$, see \cite[Thm. 2.2]{DZ_Studia2}. In other words the atoms for $H^1(\ll)$ are either classical atoms or local atoms related to some $Q_k \in \qq$.

Among examples of atoms for which one can find a family $\qq$ such that the assumptions of \cite[Thm. 2.2]{DZ_Studia2} are satisfied, there are potentials $\vv$ in the Reverse H\"older class in dimension $d\geq 3$. For more examples see \cite{DZ_Studia2}. Later, Czaja and Zienkiewicz \cite[Thm. 2.4]{CzajaZienkiewicz_ProcAMS} proved that in dimension one for any $0\leq \vv \in L^1_{loc}(\RR)$ there is a~family of intervals such that \cite[Thm. 2.2]{DZ_Studia2} gives local atomic decompositions for $H^1(\ll)$. Let us mention that results of this type in dimension $d=2$, are much unlike those in $d=1$ and $d\geq 3$, see \cite{DZ_Revista2}.

A question that we are concerned with is: what happens if we replace $-\Delta$ by the Bessel operator $\opB f(x) =  - f''(x) - \a/x \, f'(x)$ on $L^2((0,\8), x^\a\, dx)$? It is known that if $\a+1 \in \NN$ then $\opB$ corresponds to $-\Delta$ on radial functions on $\Rd$ with $d=\a+1$, however $\opB$ exists and generates a semigroup for all $\a>-1$, which can be considered as the Laplace operator on spaces with non-integer dimensions.

In this paper we prove results similar to \cite{DZ_Studia2} and \cite{CzajaZienkiewicz_ProcAMS} for the Bessel operator with a potential. We were  motivated to consider non-integer parameters $\a$ by the fact, that the Hardy space $H^1(\ll)$ admits different atomic decompositions in different dimensions, as it was mentioned above. Especially, we were interested in the dimensions $\a+1 \in (1,2)$, having in mind the difference of results for $\a=0$ and $\a=1$, see \cite{CzajaZienkiewicz_ProcAMS}, \cite{DZ_Studia2}, and \cite{DZ_Revista2}. Although we use the same scheme of proofs as in \cite{DZ_Studia2} and \cite{CzajaZienkiewicz_ProcAMS}, some technical difficulties appear. Indeed, in the space with the~weighted measure $x^\a \, dx$ the analysis is more delicate. One of the main problems is that the measure $\mu$~is not invariant under translation, and the induction argument used to prove Lemma \ref{lem:conditionD} becomes more complicated. Also, we added a precise explanation of the superharmonicity of $\phi_I$ (see \eqref{def_phi}) in the Appendix. Let us notice that the semigroup related to the Bessel operator $\opB$, which is used intensely in the proof of Theorem \ref{thm:maintheorem}, is given in terms of Bessel functions $I_\a$, but we shall use mainly the Gaussian estimates \eqref{eq:ptapprox1} for $\exp(-t\opB)$, which are obtained from asymptotics of $I_\a$. In particular, we need to prove an atomic characterization for local Hardy space related to $\opB$, see Subsection \ref{subsec22}.

Furthermore, for readers convenience, recall that with the operator $\opB$ on $L^2((0,\8), x^\a dx)$ we can relate another Bessel operator $\wt{\opB} f(x) = \Uu \opB \Uu^{-1} f(x) =  -f''(x) + \frac{\a(\a-2)}{4x^2} f(x)$, where $\Uu : L^2((0,\8), x^\a dx) \to L^2((0,\8), dx)$ is an isometry given by $\Uu f(x) = x^{\a/2} f(x)$. Therefore, the $L^2$-theory of $\opB$ and $\wt{\opB}$ can be studied simultaneously. However, it seems that $\Uu$ (or even $\Uu^2$ which is an isometry of suitable $L^1$ spaces) cannot help in studies of Hardy spaces. To see this, one can look at the atomic decompositions of the Hardy spaces $H^1(\opB)$ and $H^1(\wt{\opB})$ given in \cite{BDT}. It appears that all the atoms for $H^1(\opB)$ have cancellation conditions with respect to $x^\a dx$ (see Section \ref{sec2} below), but some atoms for $H^1(\wt{\opB})$ are of the form $a(x) =\delta^{-1} \chi_{(0,\delta)}(x)$ with $\delta >0$, so they do not satisfy any cancellations. Moreover, $H^1(\wt{\opB})$ for each $\a>0$ is the same space, which can not be said about $H^1(\opB)$.

\subsection{Definitions}\label{ssec1.1}
For $\a>0$ let $(X, \rho, \mu)$ be a metric-measure space, where $X=(0,\8)$, $\rho(x,y) = |x-y|$ and $d\mu(x) = x^\a \, dx$. Denote $B(x,r) = \set{y\in X \ : \ \rho(x,y)<r}$ and observe that $X$ is a space of homogeneous type in the sense of~Coifman-Weiss \cite{CoifmanWeiss_BullAMS}, i.e. the doubling condition holds
    \eqx{\lab{doubling}
    \mu(B(x,2r)) \leq C \mu(B(x,r)),
    }
where $C$ does not depend on $x\in X$ and $r>0$.

The classical Bessel operator is given by
\eqx{\opB f(x) =  - f''(x) - \frac{\a}{x} f'(x).}
Slightly abusing notation, we shall also write $\opB$ for the densely defined, self-adjoint operator on $L^2(X, \mu)$ that corresponds to the differential operator above, see Subsection \ref{ssec2.1} for the semigroup generated by~$\opB$.

In this paper we consider the Bessel-Schr\"odinger operator $\opL$,
\eq{
	\label{eq:operatorL} \opL f= \opB f + V\cdot f,
}
where $V \in L^1_{loc}(X,\mu), \,  V\geq 0$. To be more precise, denote $\innprod{f,g} = \int f g \, d\mu$ and define a quadratic form
\eqx{
Q(f,g) = \innprod{f', g'} + \innprod{\sqrt{V}f, \sqrt{V} g},
}
with the domain
\eq{\label{DomQ}
\Dom(Q) = \rm{cl}\set{f\in C_c^1[0,\8) \ : \ f'(0^+)=0 } \cap \set{f\in L^2(X,\mu): \sqrt{V}f\in L^2(X,\mu)},
}
where $\rm{cl}(A)$ stands for the closure of the set $A$ in the norm $\norm{f}_{L^2(X,\mu)} + \norm{f'}_{L^2(X,\mu)}$.
The quadratic form $Q$ is positive and closed. Therefore, it defines a self-adjoint operator $\opL$ with the domain
\eq{\lab{DomL}
	\Dom(\opL) = \set{f\in \Dom(Q): \exists h\in L^2(X,\mu) \ \forall g \in \Dom(Q) \quad Q(f,g)=\innprod{h,g}}.
}
For $f,h$ as above we put $\opL f:=h$. Let $\Kt{t} = \exp(-t\opL)$ be the semigroup generated by $\opL$. Denote by $B_s$ the Bessel process on $(X,\mu)$. By using the Feynman-Kac formula,
\eqx{
K_tf(x) = E^x\eee{\exp\eee{-\int_0^t V(B_s)\, ds}f(B_t)},}
one gets that $K_t$ has an integral kernel $\Kt{t}(x,y)$ and
\eq{\label{FKK}
0\leq \Kt{t}(x,y)\leq P_t(x,y),
}
where $P_t(x,y)$ is the kernel related to $P_t = \exp\eee{-t\opB}$, see Subsection \ref{ssec2.1}.

We define the Hardy space $\haL$ by means of the maximal operator associated with $\Kt{t}$, namely
\eq{
	\label{eq:hardyLdefinition} \haL = \set{f\in \Ls{X} : \norm{f}_{\haL} = \norm{\sup_{t>0} \abs{\Kt{t}f}}_{\Ls{X}} < \8}.
}
The goal of this paper is to give an atomic characterizations of local type for $\haL$. Let $|\I{}|$ be the diameter of $\I{}$.

\defn{def_proper}{
Let $\ii$ be a collection of intervals that are closed with respect to the topology on $(0,\8)$. We~call a family $\ii$ {\it a proper section} of $X$ if:
\begin{enumerate}
\item[(a)] for $I,J \in \ii$, $I\neq J$, the intersection $I\cap J$ is either the empty set or a singleton,
\item[(b)] $X=\bigcup_{I\in\ii} I$,
\item[(c)] there exists a constant $C_0>0$ such that for $I,J \in \ii$,  $I \cap J \neq \emptyset$ we have
    \eqx{C_0^{-1}|I|\leq |J| \leq C_0 |I|.}
\end{enumerate}
}
Denote $\tau B(c,r) := B(c, \tau r)$. For an interval $I = B(x,r)$ (if $I=(0,2A)$ we take $I=B(A, A)$), let $c I:= B(x, c r)$. For a family $\ii$ as in Definition \ref{def_proper} we set $\iI{}:=\be \I{}$ for some fixed $\be>1$. By choosing any $\be<\min\left(2^{1/3},(1+C_0^{-1})^{1/3}\right)$ we have that $I^{***} \cap J^{***} \neq \emptyset$ if and only if $I\cap J \neq \emptyset$ and
\eq{\label{beta}
I^{***}\subset 2I.
}

We say that a function $a: X \to \CC$ is an $(\ii,\mu)${\it -atom} if:
\begin{enumerate}
    \item[(i)] there exist $\I{} \in \ii$ and an interval $J\subset I^{**}$, such that: $\supp(a) \subset J$, $\norm{a}_\8 \leq \mu(J)^{-1}$, $\int a \, d\mu = 0$, \end{enumerate}
or
\begin{enumerate}
    \item[(ii)] there exists $\I{} \in \ii$, such that $a(x) = \mu(\I{})^{-1} \mathbbm{1}_{\I{}}(x)$.
 \end{enumerate}
The atoms as in $(ii)$ are called local atoms.

The atomic Hardy space $\hatom{}$ associated with the collection~$\ii$ is defined in the following way. We say that $f\in \hatom{}$ if
\eq{
	\label{eq:hardyatomdefinition} f(x) = \sum_{n} \la_n a_n(x),
}
where $\la_n \in \CC$, $a_n$ are $(\ii,\mu)$-atoms, and $\sum_n \abs{\la_n} <\8$.
Set
\eq{
	\label{eq:hardyatomnorm} \norm{f}_{\hatom{}} := \inf \sum_n \abs{\la_n},
}
where the infimum is taken over all possible representations of $f$ as in \eqref{eq:hardyatomdefinition}.

For a collection $\ii$ as above and $V\geq0$, $V\in L^1_{loc}(X,\mu)$ we consider the following two conditions:
\\- there exist constants $C, \e>0$ such that
\eqx{ \tag{D} \label{eq:D} \sup_{y \in \iiI{}} \int_X \Kt{2^k |I|^2}(x,y) \, d\mu(x) \leq Ck^{-1-\e} \qquad \text{for } \, I \in \ii,\, k\in\NN,}
- there exist constants $C,\de>0$ such that
\eqx{ \tag{K} \label{eq:K} \int_{0}^{2t} \int_X \Pt{s}(x,y) \mathbbm{1}_{\iiiI{}}(y) V(y) \, d\mu(y) \, ds \leq C\le\frac{t}{|I|^2}\)^{\de} \quad \text{for } \, x\in X,\, I \in \ii,\,  t\leq |I|^2.}

\subsection{Statement of results} Our first main result is the following, cf. \cite[Thm. 2.2]{DZ_Studia2}
    \thm{thm:maintheorem}{
    	Assume that a proper section $\ii$ and $0\leq V\in L^1_{loc}(X,\mu)$ are given, so that \eqref{eq:D} and \eqref{eq:K} hold. Then $\haL = \hatom{}$ and there exists a constant $C>0$, such that
    	\eqx{
    		\label{eq:hardynormseq} C^{-1} \norm{f}_{\hatom{}} \leq \norm{f}_{\haL} \leq C \norm{f}_{\hatom{}}.
    	}
}

In the second part we give an important application of Theorem \ref{thm:maintheorem}. Let us restrict ourselves to $\a\in (0,1)$. We prove that for {\bf any} $0\leq V\in L^1_{loc}(X,\mu)$ we can find a family $\ii(V)$ such that the assumptions of Theorem \ref{thm:maintheorem} hold, cf. \cite[Thm. 2.4]{CzajaZienkiewicz_ProcAMS} for the case $\a = 0$. To be more precise, let $\dd$ be a family of dyadic intervals on $(0,\8)$, that is $\mathcal{D}=\set{[k2^n,(k+1)2^n]:k\in\NN, n\in\ZZ}\cup\set{(0,2^n]: n\in \ZZ}$. Consider the family $\ii(V)$ that consists of maximal dyadic closed intervals $I$ that satisfy
	\eqx{ \tag{S} \label{eq:S} \frac{\abs{2\I{}}^2}{\mu(2\I{})} \int_{2\I{}} V(y) \, d\mu(y) \leq 1.}

In Section \ref{sec4} we prove that $\ii(V)$ is a well defined proper section. The second main result is the following.
\thm{thm:maintheorem2}{}
Let $\a \in (0,1)$ and $0\leq V\in L^1_{loc}(X,\mu)$. Then the family $\ii(V)$ satisfies the assumptions of Theorem \ref{thm:maintheorem}.

\cor{coro_main}{
Let $\a \in (0,1)$ and $0\leq V\in L^1_{loc}(X,\mu)$. Then there is $C>0$, such that
    \eqx{
    C^{-1}\norm{f}_{\hh^1_{at}(\ii(V), \mu)} \leq \norm{f}_{\haL} \leq C \norm{f}_{\hh^1_{at}(\ii(V), \mu)}.}
}
The paper is organized as follows. In Section \ref{sec2} we study the atomic Hardy spaces related to $\opB$ and its local versions. This is used in a proof of Theorem \ref{thm:maintheorem}, which is provided in Section \ref{sec3}. Finally, in Section \ref{sec4} a proof of Theorem \ref{thm:maintheorem2} is given.

\section{Hardy spaces for the Bessel operator}\label{sec2}
\subsection{Global Hardy space for $\opB$}\lab{ssec2.1}
In this section we consider the case $\a>0$. Let $\Pt{t}=\exp(-t\opB)$ be the Bessel semigroup given by
\al{\nonumber
    &\Pt{t}f(x) = \int_X \Pt{t}(x,y) f(y)\,d\mu(y),\\
	&\label{eq:besselkernel} \Pt{t}(x,y) = (2t)^{-1}\exp\eee{-\frac{x^2+y^2}{4t}}I_{\ap}\le\frac{xy}{2t}\)(xy)^{-\ap},
	}
where $I_\a(x)=\sum_{m=0}^{\8} \frac{1}{m!\Gamma(m+\a+1)}\le\frac{x}{2}\)^{2m+\a}$ is the modified Bessel function of the first kind. It is clear that $\Pt{t}(x,y) = \Pt{t}(y,x)$ and, since $\opB \mathbbm{1}_{(0,\8)}(x) = 0$, we have that
    \eq{\label{contraction}
    \int_X \Pt{t}(x,y) \, d\mu(x) = 1.
    }
Let us recall that
    \eq{\lab{ballmeasure}
    \mu(B(x,t)) \simeq t(x+t)^\a.
    }
It is known that the kernel $P_t(x,y)$ satisfies the two-side Gaussian estimates (see, e.g. \cite[Lem. 4.2]{DPW_JFAA}),
    \eq{\label{eq:ptapprox1}
    C^{-1} \mu(B(x,\sqrt{t}))^{-1} \exp\eee{-\frac{\abs{x-y}^2}{c_1t}} \leq \Pt{t}\xy \leq C \mu(B(x,\sqrt{t}))^{-1} \exp\eee{-\frac{\abs{x-y}^2}{c_2t}},
    }
while the derivative satisfies
    \eq{\label{gauss_prim}
    \abs{\frac{\d}{\d x} \Pt{t}\xy} \leq C t^{-1/2} \mu(B(x,\sqrt{t}))^{-1} \exp\eee{-\frac{\abs{x-y}^2}{ct}}.
    }

Let $H^1(\opB)$ be the Hardy space related to $P_t$, i.e. the space defined as in \eqref{eq:hardyLdefinition} with $\opL$ and $K_t$ replaced by $\opB$ and $P_t$, respectively. We call a function $a$ an $\mu$-{\it atom} if
\begin{enumerate}
\item[(iii)] there exists an interval $J\subseteq X$, such that: $\supp(a) \subset J$, $\norm{a}_\8 \leq \mu(J)^{-1}$,  $\int a \, d\mu = 0$.
\end{enumerate}
The atomic Hardy space $H^1_{at}(\mu)$ is defined as in \eqref{eq:hardyatomdefinition} and \eqref{eq:hardyatomnorm} with $a_n$ being $\mu$-atoms.

\thm{thm:hardybessel}{
\cite[Thm. 1.7]{BDT} Let $\a >0$. There is $C>0$ such that
    \eqx{ C^{-1}\norm{f}_{H^1_{at}(\mu)} \leq \norm{f}_{H^1(\opB)} \leq C \norm{f}_{H^1_{at}(\mu)}.}
}

\subsection{Local Hardy space for $\opB$}\label{subsec22}
For $\tau >0$ we define $h^1_\tau(\opB)$, {\it the local Hardy space related to $\opB$}, as the set of $L^1(X,\mu)$ functions for which the~norm
\eqx{
	\label{eq:hardylocalnorm} \norm{f}_{\ha{\tau}} = \norm{\st{\tau} \abs{\Pt{t}f}}_{\Ls{X}}
}
is finite.

Let $\ii_{\tau}$ be a proper section of $X$ that consists of closed intervals of length $\tau$.

\thm{thm:atomiclocalnormseq}{
    $(a)$ There exists $C>0$ such that for $\tau>0$ and an $(\ii_{\tau},\mu)$-atom $a$ we have
        \eqx{\norm{a}_{\ha{\tau}} \leq C.}

    $(b)$ Let $I$ be an interval such that $\supp(f) \subset I^{*}$ and $f\in\ha{|I|}$. Then
	\eqx{ f = \sum_{n=0}^\8 \la_n a_n, \quad \sum_{n=0}^\8 \abs{\la_n} \leq C \norm{f}_{\ha{|I|}},}
	where $a_0(x) = \mu(I)^{-1} \mathbbm{1}_I(x)$ and $a_n$ are $\mu$-atoms supported in $I^{**}$ for $n\geq 1$.
}
Let us remark that another characterization of $h^1_{\tau}(\opB)$, by mean of a local Riesz transform, was given in \cite[Thm. 2.11]{Preisner_JAT}. Theorem \ref{thm:atomiclocalnormseq} will be used to prove Theorem \ref{thm:maintheorem}.

\pr{
$(a)$ Obviously, if $a$ is $\mu$-atom, then the statement follows from Theorem \ref{thm:hardybessel}. Assume then that $a(x) = \mu(I)^{-1} \mathbbm{1}_I(x)$ and $|I|=\tau$. It is well known that \eqref{eq:ptapprox1} implies the boundedness of the maximal operator $\sup_{t>0} \abs{\Pt{t} f}$ on $L^2(X,\mu)$. Using this fact and the Schwarz inequality,
    \eqx{
    \norm{\sup_{t\leq \tau^2} \abs{\Pt{t} a}}_{L^1(I^{**},\mu)} \leq C \mu(I)^{1/2} \norm{\sup_{t>0} \abs{\Pt{t} a}}_{L^2(X,\mu)} \leq C \mu(I)^{1/2} \norm{a}_{L^2(X,\mu)} \leq C.
    }
Denote by $c_I$ the center of $I$ and notice that $|x-y| \simeq |x-c_I|$ when $y\in I$ and $x\in (I^{**})^c$. By \eqref{eq:ptapprox1},
    \alx{
    \norm{\sup_{t\leq \tau^2} \abs{\Pt{t} a}}_{L^1((I^{**})^c,\mu)}
    &\leq C \int_{(I^{**})^c} \sup_{t\leq \tau^2} \int_{I} \mu(B(x,\sqrt{t}))^{-1} \exp\eee{-\frac{|x-c_I|^2}{c_3 t}} \mu(I)^{-1} \, d\mu(y) \, d\mu(x)\\
    &\leq C \int_{(I^{**})^c} \sup_{t\leq |I|^2} t^{-1/2} \frac{t}{|x-c_I|^2} \, dx\\
    &\leq C \int_{(I^{**})^c} |I| |x-c_I|^{-2} \, dx \leq C.
    }
$(b)$ Define $\la_0 := \int f d\mu$ and $g(x) = f(x) - \la_0 \mu(I)^{-1} \mathbbm{1}_I(x)$. Notice that
    \eqx{
    \abs{\la_0} \leq \norm{f}_{L^1(X,\mu)} \leq \norm{f}_{\ha{|I|}}.
    }
Therefore, $\norm{g}_{L^1(X,\mu)} \leq 2 \norm{f}_{L^1(X,\mu)}$. Our goal is to prove that
    \eq{\label{plk}
    \norm{g}_{H^1(\opB)} \leq C \norm{f}_{\ha{|I|}}.
    }
For $t\leq |I|^2$ using $(a)$ we obtain
    \eqx{
    \norm{\sup_{t\leq |I|^2} \abs{\Pt{t} g}}_{L^1(X,\mu)} \leq \norm{f}_{\ha{|I|}} +|\la_0|\norm{\mu(I)^{-1}\mathbbm{1}_I}_{\ha{|I|}} \leq C \norm{f}_{\ha{|I|}}.
    }
For $t\geq |I|^2$ we shall use $\int g d\mu = 0$. By the symmetry of $P_t(x,y)$ and \eqref{gauss_prim},
    \alx{
    \norm{\sup_{t\geq |I|^2} \abs{\Pt{t} g}}_{L^1((I^{**})^c,\mu)} &\leq \int_{(I^{**})^c} \sup_{t\geq |I|^2} \abs{\int_{I^*}\eee{\Pt{t}(x,y) - \Pt{t}(x,c_I)}g(y) \, d\mu(y)} \, d\mu(x)\\
    &\leq \int_{(I^{**})^c} \int_{I^*} \sup_{t\geq |I|^2} \frac{|y-c_I|}{\sqrt{t}}\mu(B(x,\sqrt{t}))^{-1} \exp\eee{-\frac{|x-c_I|^2}{c_2 t}}\abs{ g(y)} \, d\mu(y) \, d\mu(x)\\
    &\leq C \norm{g}_{L^1(X,\mu)} |I| \int_{(I^{**})^c} \sup_{t\geq |I|^2} t^{-1} \exp\eee{-\frac{|x-c_I|^2}{c_2 t}} \, dx\\
    &\leq C \norm{f}_{L^1(X,\mu)} |I| \int_{(I^{**})^c} |x-c_I|^{-2} \, dx \leq C \norm{f}_{\ha{|I|}}.
    }
Likewise,
    \alx{
    \norm{\sup_{t\geq |I|^2} \abs{\Pt{t} g}}_{L^1(I^{**},\mu)} &\leq C \int_{I^{**}} \int_{I^*} \sup_{t\geq |I|^2} \mu(B(x,\sqrt{t}))^{-1} \abs{g(y)}\, d\mu(y) \, d\mu(x)\\
    &\leq C \norm{g}_{L^1(X,\mu)} \int_{I^{**}} \sup_{t\geq |I|^2} t^{-1/2} \, dx \leq C  \norm{f}_{\ha{|I|}}.
    }

From \eqref{plk} we have that $g\in H^1(\opB)$, so using Theorem \ref{thm:hardybessel} we obtain $\la_k$ and $\mu$-atoms $a_k$ such that \break $g=\sum_{k=1}^\8 \la_k a_k$. Consequently,
        \eqx{
        f=\sum_{k=0}^{\8} \la_k a_k, \qquad \sum_{k=0}^\8 \abs{\la_k} \leq C\norm{f}_{h^1_{|I|}(\opB)},
        }
        where $a_0(x)=\mu(I)^{-1}\mathbbm{1}_I(x)$.

The only problem we have to deal with is that $a_k$ are not necessarily supported in $I^{**}$. Let $\psi$ be a function such that $\psi \equiv 1$ on $I^{*}$, $\psi \equiv 0$ on $(I^{**})^c$, and $\norm{\psi'}_\8 \leq C|I|^{-1}$. Then $f=\psi f$. To complete the proof we will show, that for every $\mu$-atom $a_k$ there exist sequences $\la_k^j\in\CC$  and $a_k^{j}$, such that each $a_k^j$ is either $\mu$-atom or $a_k^j =\mu(I)^{-1}\mathbbm{1}_I$, and $ a_k = \sum_j \la_k^{j} a_k^{j}$, $\sum_j |\la_k^{j}| < C$, with $C$ that not depend on $k$. Fix $a=a_k$ and an interval $J\subseteq X$ such that: $\supp(a) \subset J$, $\norm{a}_\8 \leq \mu(J)^{-1}$. Obviously $\supp(\psi a) \subseteq I^{**} \cap J$, and if $J\subset I^{*}$, then $\psi a = a$  is an $\mu$-atom. Furthermore, if $J\subset (I^{**})^c$, then $\psi a = 0$, so it suffices to consider the case that $J \cap (I^{*})^c\cap I^{**} \neq \emptyset$. Denote $K:=I^{**}\cap J$. Observe that $|K|\leq C_0|I|$ and let $N\in\NN$ be such that $1/2^{N+1}|I| \leq C_0^{-1}|K| \leq 1/2^{N}|I|$.

Define $\la := \int \psi a \, d\mu$  and notice that
	\alx{ |\la| = \abs{\int_K \psi(x) a(x) \, d\mu(x)} & = \abs{\int_K a(x)(\psi(x) - \psi(c_{K})) \, d\mu(x)} \\
	& \leq C |I|^{-1} \mu(K)^{-1} \int_K |x-c_{K}| \, d\mu(x) \\
	& \leq C |I|^{-1} |K|  \leq C 2^{-N}. }
Let us choose intervals $I_j$, such that $K =: I_0 \subset I_1 \subset \cdots \subset I_{N} \subset I^{**}$, where $|I_{j+1}|/|I_j|=2$ and $|I_N|\simeq |I|$. Then
    \alx{ \psi a & = \left( \psi a - \la \mu(I_0)^{-1}\mathbbm{1}_{I_0} \right) + \sum_{j=1}^{N} \la \left( \mu(I_{j-1})^{-1}\mathbbm{1}_{I_{j-1}} - \mu(I_{j})^{-1}\mathbbm{1}_{I_{j}} \right) \\
    & + \la (\mu(I_N)^{-1}\mathbbm{1}_{I_N} - \mu(I)^{-1}\mathbbm{1}_{I}) + \la \mu(I)^{-1}\mathbbm{1}_I = \sum_{j=0}^{N+2} b_j.\\ \quad &
    }
    Observe that:
    \begin{enumerate}
    \item[(1)] $\supp(b_0)\subseteq I_0$, $\int b_0 \, d\mu = 0$, and
    \eqx{ \norm{b_0}_\8 \leq C(1+2^{-N})\mu(I_0)^{-1} \leq C \mu(I_0)^{-1},}
	\item[(2)] for $j=1,...,N$ we have: $\supp(b_j) \subseteq I_j$, $\int b_j \, d\mu = 0$, and
    \eqx{\norm{b_j}_\8 \leq C |\la| \mu(I_j)^{-1} \leq C 2^{-N} \mu(I_j)^{-1},}
	\item[(3)] $\supp(b_{N+1}) \subseteq I^{**}$, $\norm{b_{N+1}}_\8 \leq C |\la| \mu(I)^{-1}$ and $\int b_{N+1} \, d\mu = 0$.
	\end{enumerate}
We conclude that $b_j$ are multiples of $(\ii_{|I|},\mu)$-atoms and
	\alx{
	\norm{\sum_{j=0}^{N+2} b_j}_{H^1_{at}(\ii_{|I|},\mu)} \leq C \sum_{j=0}^{N+2} 2^{-N} \leq C.
	}
}

\cor{cororor}{
There exists a constant $C>0$ such that for $\tau>0$ we have
	\eqx{ C^{-1}\norm{f}_{\hatom{\tau}} \leq \norm{f}_{h^1_\tau(\opB)} \leq C \norm{f}_{\hatom{\tau}}.}
}
The right inequality in Corolarry \ref{cororor} follows easily from Theorem \ref{thm:atomiclocalnormseq}$(a)$. For the left inequality one uses Theorem \ref{thm:atomiclocalnormseq}$(b)$ with a suitable partition of unity and methods as in Lemmas \ref{lem:intfunctionwithPt} and \ref{lem:diffkt} below. We omit the details.

\section{Proof of Theorem \ref{thm:maintheorem}}\label{sec3}

\subsection{Auxiliary estimates}

For a proper section $\ii$ let $\set{\p{I}}_{I\in \ii}$ be a partition of unity associated with $\ii$, that is a family of $C^\8$ functions on~$X$, such that $\supp(\p{I})\subset I^{*},$
$ 0 \leq \p{I} \leq 1,\, \norm{\p{I}'}_\8 \leq |I|^{-1} $ and $\sum_{I\in\ii}\p{I}(x) = 1$ for all $x \in X$.

The perturbation formula states that
    \eq{\lab{perturbation}
    \Pt{t}\xy - \Kt{t}\xy = \int_0^t \int_X \Pt{t-s}(x,z) V(z) \Kt{s}(z,y) \, d\mu(z) ds.
    }

To prove Theorem \ref{thm:maintheorem} we closely follow the proof of \cite[Thm. 2.2]{DZ_Studia2}. However, in our weighted space, some technical difficulties appear. Therefore we present all the details for convenience of the reader.

\lem{lem:intfunctionwithPt}{ For $I\in\ii$ and $f\in\Ls{X}$,
\eqx{ \norm{\st{4|I|}\abs{\Pt{t}(\p{I}f)}}_{\Ls{(\iiI{})^c}} \leq C \norm{\p{I}f}_{\Ls{X}}.}
}

\pr{
	Denote by $c_{I}$ the center of the interval $\I{}\in\ii$. For $x\in(\iiI{})^c$ and $y\in\iI{}$ we have $\abs{x-y} \simeq \abs{x-c_I}$. Notice that \eqref{ballmeasure} implies
    \eq{\label{uuu}
    \mu(B(x,\sqrt{t}))^{-1} x^\a \leq t^{-1/2}.
    }
Using \eqref{eq:ptapprox1} the left-hand side is bounded by
	\spx{
\int_{(\iiI{})^c} \st{4|I|}  &  \int_{I^*} \mu(B(x,\sqrt{t}))^{-1} e^{-\frac{\abs{x-y}^2}{c_2t}} \abs{\p{I}(y)f(y)} \, d\mu(y) \, d\mu(x) \\
		& \leq C \norm{\p{I} f}_{L^1(X,\mu)} \int_{(\iiI{})^c} \st{4|I|} x^{-\a} t^{-\frac{1}{2}} e^{-\frac{\abs{x-c_I}^2}{c_3 t}} \, d\mu(x) \\
		& \leq C \norm{\p{I} f}_{L^1(X,\mu)} \int_{(\iiI{})^c}  \st{4|I|}t^{-\frac{1}{2}} \frac{t}{|x-c_I|^2} \, dx \\
		& \leq C \norm{\p{I} f}_{L^1(X,\mu)} \int_{(\iiI{})^c} |I| |x-c_I|^{-2} \, dx \\
		& \leq C \norm{\p{I}f}_{\Ls{X}}.
	}
}

\lem{lem:firstlemma}{
	For $I\in\ii$ and $f\in\Ls{X}$,
	\eqx{ \norm{\p{I} f}_{h^1_{|I|}(\opB)} \simeq \norm{ \st{|I|} |\Pt{t}(\p{I} f)|}_{\Ls{I^{**}}}. }
	}

\newcommand{\ov}[1]{\bar{#1}}
Denote
    \eqx{
    \wt{f}_I:=\sum_{J:\I{}\cap J \neq \emptyset} \p{J} f, \qquad \qquad \ov{f}_I := f-\wt{f}_I = \sum_{J:\I{}\cap J = \emptyset} \p{J} f.
    }

\lem{lem:diffkt}{
	For $I\in\ii$ and $f\in\Ls{X}$,
    \eqx{
    \norm{\st{|I|}\abs{{\Kt{t}} \left(\p{I} \wt{f}_I\right) - \p{I} \Kt{t} \left(\wt{f}_I\right)}}_{\Ls{\iiI{}}}  \leq C \sum_{J:\I{}\cap J \neq \emptyset} \norm{\p{J} f}_{\Ls{X}}.
    }
}

\pr{
	Denote $\wt{I} = \cup_{J:I\cap J \neq \emptyset} J^{*}$. Note that for $x\in\iiI{}$ and $y\in \wt{I}$, there is  $\abs{x-y}\leq C |I|$ and
    \eq{\label{eq:jjj}
    \st{|I|}~t^{-\frac{1}{2}}~e^{-\frac{\abs{x-y}^2}{c_2t}}~\leq~C~\abs{x-y}^{-1}.
    }
    Using \eqref{FKK},  \eqref{eq:ptapprox1}, \eqref{uuu}, and \eqref{eq:jjj},

	\spx{
		\norm{\st{|I|}\abs{\Kt{t}\left(\p{I}\wt{f}_I\right)-\p{I}\Kt{t}\left(\wt{f}_I\right)}}_{\Ls{\iiI{}}} & \leq \int_{\iiI{}} \st{|I|} \int_{\wt{I}} \abs{\p{I}(y) -\p{I}(x)} \Kt{t}(x,y) \abs{\wt{f}_I(y)} \, d\mu(y) \, d\mu(x) \\
		& \leq C \int_{\iiI{}} \int_{\wt{I}} \abs{x-y}\norm{\p{I}'}_\8  \st{|I|} t^{-\frac{1}{2}}  e^{-\frac{\abs{x-y}^2}{c_2 t}} \abs{\wt{f}_I(y)} \, d\mu(y) \, dx \\
		& \leq C \int_{I^{**}} |I|^{-1} dx \norm{\wt{f}_I}_{L^1(X,\mu)}.
	}
}

\lem{lem:suppthreestars}{
	Assume that $V$ and $\ii$ are given, so that \eqref{eq:D} holds. Then
	\eqx{\sum_{I\in\ii} \norm{\sup_{t>0} \abs{\Kt{t} \left( \ov{f}_I\right) }}_{\Ls{\iiiI{}}} \leq C \norm{f}_{\Ls{X}}.}
}

\pr{
	Denote $s_m=2^m |I|^2$ and let $m\geq 2$. By the semigroup property, \eqref{FKK}, \eqref{eq:ptapprox1},
	\spx{
		&\norm{\sup_{s_m \leq t \leq s_{m+1}}  \Kt{t}(\abs{\p{J}f})}_{\Ls{X}} \\
		& \leq \int_{X} \sup_{s_m\leq t \leq s_{m+1}} \int_X \int_X \Pt{t-s_{m-1}}(x,z) \Kt{s_{m-1}}(z,y) \abs{\p{J}(y)f(y)} \, d\mu(z) \, d\mu(y) \, d\mu(x) \\
		& \leq C   \int_X \abs{\p{J}(y)f(y)} \int_X \Kt{s_{m-1}}(z,y)  \int_{X} s_m^{-1/2} \exp\expr{-\frac{|x-z|^2}{c_3 s_m}} dx \, d\mu(z) \, d\mu(y)\\
		& \leq C (m-1)^{-1-\e} \cdot \norm{\p{J}f}_{\Ls{X}}.
	}
	In the last inequality we have applied \eqref{eq:D}. Using the above estimate and Lemma \ref{lem:intfunctionwithPt},
	\spx{
		&\sum_{I\in\ii} \norm{\sup_{t>0} \abs{\Kt{t} \left( \sum_{J:I\cap J=\emptyset} \p{J} f \) }}_{\Ls{\iiiI{}}} \\	
		& \leq \sum_{J\in\ii} \sum_{I:J\cap I = \emptyset} \norm{\sup_{t>0}\Kt{t}(\abs{\p{J}f})}_{\Ls{\iiiI{}}} \leq C \sum_{J\in\ii} \norm{\sup_{t>0}\Kt{t}(\abs{\p{J}f})}_{\Ls{(J^{**})^c}} \\
		& \leq C \sum_{J\in\ii} \norm{\st{4|J|}\Kt{t}(\abs{\p{J}f})}_{\Ls{(J^{**})^c}} + C  \sum_{J\in\ii} \sum_{m=2}^\8 \norm{\sup_{s_m\leq t \leq s_{m+1}}\Kt{t}(\abs{\p{J}f})}_{\Ls{(J^{**})^c}} \\
		& \leq C \sum_{J\in\ii} \norm{\p{J}f}_{\Ls{X}} + C  \sum_{J\in\ii} \sum_{m=2}^\8 m^{-1-\e} \cdot \norm{\p{J}f}_{\Ls{X}}  \leq C \norm{f}_{\Ls{X}}.
	}
}

\lem{lem:pertur}{
	\eqx{\int_0^\8 \int_X   V(z)\Kt{s}(\abs{f})(z)  \, d\mu(z)\,ds \leq C \norm{f}_{\Ls{X}} }
}
\pr{
	By integrating \eqref{perturbation} and using \eqref{contraction} we obtain
	\eqx{  \int_0^t \int_X  V(z) \Kt{s}(z,y) \, d\mu(z) \, ds \leq C.}
    The proof is finished by setting $t\to \8$.
}

\lem{lem:diffptandkt}{
	 Assume that $V$ and $\ii$ are given, so that \eqref{eq:K} holds. Then for $I\in\ii$,
	\eqx{
	 \norm{\st{|I|} \abs{(\Pt{t}-\Kt{t})(\p{I} f)}}_{\Ls{X}} \leq C \norm{\p{I} f}_{\Ls{X}}.
	}
}

\pr{
	By \eqref{FKK} and Lemma \ref{lem:intfunctionwithPt},
	\eqx{  \norm{\st{|I|} \abs{(\Pt{t}-\Kt{t})(\p{I} f)}}_{\Ls{(\iiI{})^c}} \leq C \norm{\p{I} f}_{\Ls{X}}.}
Consider the integral on $I^{**}$. From \eqref{perturbation},
	\eqx{(\Pt{t} - \Kt{t})(\p{I}f)(x) = \int_0^t\Pt{t-s}(V''\Kt{s}(\p{I}f))(x) \, ds + \int_0^t\Pt{t-s}(V'\Kt{s}(\p{I}f))(x) \, ds }
	where $V=\mathbbm{1}_{\iiiI{}}V + \mathbbm{1}_{(\iiiI{})^c}V = V' + V''$. Repeating the argument from the proof of Lemma \ref{lem:intfunctionwithPt} and using Lemma \ref{lem:pertur} we obtain
	\eqx{\norm{\st{|I|} \int_0^t\Pt{t-s}(V''\Kt{s}(\p{I}f)) \, ds }_{\Ls{\iiI}} \leq C \norm{\p{I}f}_{\Ls{X}}.}
	To estimate the integral that contains $V'$ write
	\eqx{ \abs{\int_0^t\Pt{t-s}(V'\Kt{s}(\p{j}f))(x) \, ds}  \leq \int_0^{t/2}\Pt{t-s}(V'\Pt{s}(\abs{\p{j}f}))(x) \, ds + \int_{t/2}^t\Pt{t-s}(V'\Pt{s}(\abs{\p{j}f}))(x) \, ds}
	\eqx{ = A_{1,t}(x) + A_{2,t}(x).}

	Let $t_m = 2^{-m} |I|^2$. Similarly as in proof of the Lemma \ref{lem:suppthreestars} we obtain	
	\spx{
		\norm{ \st{|I|} A_{1,t}}_{\Ls{X}} & \leq \sum_{m=0}^\8 \int_X \sup_{t_{m+1} \leq t \leq t_m} \int_0^{t/2} \int_X \Pt{t-s}(x,y)V'(y)\Pt{s}(\abs{\p{I}f})(y) \, d\mu(y) \, ds \, d\mu(x) \\
	    & \leq C \sum_{m=0}^\8 \int_X  \int_0^{t_m}   V'(y)\Pt{s}(\abs{\p{I}f})(y)  \int_X t_m^{-\frac{1}{2}}  e^{-\frac{\abs{x-y}^2}{{c_3 t_m}}}  \, dx \, ds\,  d\mu(y) \\
	    & \leq C \sum_{m=0}^\8 \int_0^{t_m} \int_X  \mathbbm{1}_{\iiiI{}}(y)V(y)\Pt{s}(\abs{\p{I}f})(y) \,  d\mu(y)  \, ds \\
	    & \leq C  \sum_{m=0}^\8 \left( \frac{2^{-m}|I|^2}{|I|^2} \right) ^{\de} \norm{\p{I}f}_{\Ls{X}}.
    }
    In the last inequality we have used \eqref{eq:K}. To estimate $A_{2,t}$ we proceed similarly noticing that for $t\in [t_{m+1}, t_m]$ and $s\in [t/2,t]$ we have $s \in [t_{m+2},t_m]$. The details are left to the reader.
}

\subsection{Proof of Theorem \ref{thm:maintheorem}}
	
\textbf{First inequality.} Assume that $f\in\haL$. Observe that $\p{I} f = \p{I}\wt{f}_I$ and
    \eqx{
    \Kt{t} \expr{\p{I} f} = \Kt{t} \expr{\p{I} \wt{f}_I} - \p{I} \Kt{t}\eee{\wt{f}_I} - \p{I} \Kt{t}\eee{\ov{f}_I}+ \p{I} \Kt{t}\eee{f}.
    }
From Lemmas \ref{lem:diffkt} and \ref{lem:suppthreestars} we deduce that

	\spx{
		 \sum_{I\in\ii} \norm{\st{|I|} \abs{\Kt{t}(\p{I} f)}}_{\Ls{\iiI{}}}
		\leq & \sum_{I\in\ii} \norm{\st{|I|} \abs{\Kt{t}\left( \p{I} \wt{f}_I \right) - \p{I} \Kt{t}  \left( \wt{f}_I \right) }}_{\Ls{\iiI{}}} \\
		& + \sum_{I\in\ii}\norm{ \st{|I|} \abs{ \p{I} \Kt{t}\left( \ov{f}_I \)}}_{\Ls{\iiI{}}} \\
		& + \sum_{I\in\ii} \norm{ \st{|I|} \abs{ \p{I} \Kt{t}f }}_{\Ls{\iiI{}}} \\
    \leq& C  \eee{\norm{f}_{\Ls{X}}+ \norm{f}_{\haL}} \leq C\norm{f}_{\haL}.
	}

The above estimate, together with Corollary \ref{lem:firstlemma} and Lemma \ref{lem:diffptandkt}, lead to
	
	\spx{
			\sum_{I\in\ii} \norm{\p{I} f}_{\ha{|I|}} 
			& \leq C \sum_{I\in\ii} \norm{\st{|I|}\abs{(\Pt{t}-\Kt{t})(\p{I} f)}}_{\Ls{\iiI{}}} + C\sum_{I\in\ii} \norm{\st{|I|} \abs{\Kt{t}(\p{I} f)}}_{\Ls{\iiI{}}} \\
			& \leq C \norm{f}_{\haL}. \\
		}
Now we use Theorem \ref{thm:atomiclocalnormseq}$(b)$ for each $\p{I} f$ getting $(\la_n^I)_n$ and $(\ii, \mu)$-atoms $(a_n^I)_n$, so that
    \eqx{\p{I}f(x)=\sum_{n} \la_n^I a_n^I(x), \quad \text{and} \quad  \sum_{n} \abs{\la_n^I} \leq C \norm{\p{I} f}_{\ha{|I|}}.}
Notice that $\sum_{n,I}\norm{\la_n^Ia_n^I}_{\Ls{X}} < \8$, since $\norm{a_n^I}_{\Ls{X}}\leq 1$. Summing up for all $I\in \ii$ we finish the first part of the proof.


	\textbf{Second inequality.} Let $a$ be an $(\ii,\mu)$-{\it atom}, such that $\supp(a) \subset \iiI{}$. There exists an integer $ m \geq 0$, independent of $I$, such that
		\eqx{\inf \set{|I|^2 : J\in\ii, \, J\cap I \neq \emptyset} \geq 2^{-m}|I|^2.}
	Denote $t_n=2^{-n}|I|^2$, $n\in\ZZ$.
	Observe that
	\eqx{\norm{\sup_{t\leq t_m} \abs{\Kt{t}(a)}}_{\Ls{X}} \leq \norm{\sup_{t\leq t_m} \abs{(\Kt{t}-\Pt{t})(a)}}_{\Ls{X}} + \norm{\sup_{t\leq t_m} \abs{\Pt{t}( a)}}_{\Ls{X}}\leq \norm{a}_{\Ls{X}}.}
	In the last inequality we applied Lemma \ref{lem:diffptandkt} and Theorem \ref{thm:atomiclocalnormseq}$(b)$, since $a$ is also an $(\ii_{|I|},\mu)$-{\it atom}.

	It suffices to estimate $\norm{\sup_{t>t_m}\abs{\Kt{t}a}}_{\Ls{X}}$. This is done by using \eqref{FKK} and \eqref{eq:K}. Indeed, using similar  methods as in the proof of Lemma \ref{lem:intfunctionwithPt},
	\spx{
		\norm{\sup_{t>t_m} \abs{\Kt{t}a}}_{\Ls{X}} & \leq \sum_{n\leq m} \norm{\sup_{t_{n} \leq t \leq t_{n-1}}\abs{\Kt{t}a}}_{\Ls{X}} \\ & \leq \sum_{n\leq m} \norm{\sup_{t_{n+1}\leq t \leq 3t_{n+1}}\Kt{t}\abs{\Kt{t_{n+1}}a}}_{\Ls{X}} \\
		& \leq C \sum_{n\leq m+1} \int_X \int_X \Kt{2^{-n}|I|^2}(x,y)\abs{a(y)} \, d\mu(y)  \, d\mu(x) \\
		& \leq C \norm{a}_{\Ls{X}}\leq C,
	}
since the operator $\sup_{t_{n+1}\leq t \leq 3t_{n+1}} K_t$ is bounded on $L^1(X,\mu)$, see \eqref{eq:ptapprox1}.

\section{Proof of Theorem \ref{thm:maintheorem2}.}\label{sec4}

In the whole section we assume that $\a\in(0,1)$. Recall that for $0\leq a<b$,
    \eq{\lab{asympt_2}
    \mu((a,b)) = \frac{ b^{1+\a} - a^{1+\a}}{1+\a} \simeq
    \begin{cases}
    b^{\a+1} & 2a\leq b\\
    (b-a) a^{\a} & 2a\geq b.
    \end{cases}
    }
If $I=(a,b)\subset(0,\8)$, by \eqref{asympt_2} and the Mean-Value Theorem we easily get
    \eq{\label{costam}
    \frac{|I|^2}{\mu(I)} \simeq b^{1-\a} - a^{1-\a}.
    }

\lem{enlargement}{
For $\a\in(0,1)$ and $I\subset J \subset (0,\8)$,
\eqx{\frac{\abs{I}^2}{\mu(I)} \int_{I} V(y) \, d\mu(y) \leq \frac{\abs{J}^2}{\mu(J)} \int_{J} V(y) \, d\mu(y).}
}

\pr{
Obviously, since $V\geq 0$, it is enough to prove that
    \eq{\lab{intinint}
    \frac{\abs{I}^2}{\mu(I)} \leq \frac{\abs{J}^2}{\mu(J)},
    }
provided that $I\subset J \subset (0,\8)$.

Let $a\leq b\leq c \leq d$ and $I=(b,c)\subset (a,d) =J$. Denote $\Gamma(x,y) = (y-x)^2/(y^{\a+1} - x^{\a+1})$. Now, \eqref{intinint} is equivalent to $\Gamma(b,c)\leq \Gamma(a,d)$. This is done in two steps.

{\bf Step 1:} $\Gamma(b,c) \leq \Gamma(a,c)$. Denote $a=sc$ and $b=tc$, where $0<s\leq t<1$. It is enough to prove
    \eqx{
    \frac{(1-t)^2}{1-t^{\a+1}} \leq \frac{(1-s)^2}{1-s^{\a+1}}.
    }
By a simple calculus argument, the function $F_1(x) = (1-x)^2/(1-x^{\a+1})$ is monotonically decreasing for $x\in (0,1)$.

{\bf Step 2:} $\Gamma(a,c) \leq \Gamma(a,d)$. Similarly, let $c=t'a$, $d=s'a$, $1<t'\leq s'$. The function \break $F_2(x) = (x-1)^2/(x^{\a+1}-1)$ is monotonically increasing in $(1,\8)$, thus
    \eqx{
    \frac{(t'-1)^2}{(t')^{\a+1}-1} \leq \frac{(s'-1)^2}{(s')^{\a+1}-1}.
    }
}

\prop{construction}{
Let $V\in L^1_{loc}(X,\mu)$, $V\geq 0$. Then the family $\ii(V)$ of maximal dyadic intervals satisfying \eqref{eq:S} is a proper section (see Definition \ref{def_proper}).
}

\pr{
For a closed dyadic interval $I$ consider $F(I)= |2I|^2 \mu(2I)^{-1} \int_{2I} V d\mu$ and denote by $\dI{}$ the smallest dyadic interval properly containing $\I{}$. Notice that $2I \subseteq 2\dI{}$ and, by Lemma \ref{enlargement}, we have $F(I)\leq F(\dI{})$. Also, for an increasing sequence of dyadic intervals $I_n \subset I_{n+1}$ we have $\lim_{n\to \8} F(I_n) = \8$ and $\lim_{n\to -\8} F(I_n) = 0$, see \eqref{costam}.

This justifies the choice of $\ii(V)$ as maximal dyadic intervals such that \eqref{eq:S} holds. What is left to prove is that $\ii(V)$ is a proper section, namely we need to show that for $I,J \in \ii(V)$, $I\cap J \neq \emptyset$ we have $|I| \simeq |J|$.

By contradiction, suppose that there exist $I_k$, $J_k$ such that $I_k \cap J_k \neq \emptyset$ and $|I_k|/|J_k| \to \8$. We can assume that $2J_k^d \subseteq 2 I_k$ for all $k$. Denote, $a_k = |I_k|^2 \mu(I_k)^{-1} |J_k|^{-2} \mu(J_k)$. By the choice of $\ii$,
	\eqx{1\geq \frac{\abs{2I_k}^2}{\mu(2\I{k})} \int_{2\I{k}} V(y) \, d\mu(y) \geq \frac{\abs{2I_k}^2\mu(2J_k^d)	}{\mu(2I_k)\abs{2J_k^d}^2} \frac{\abs{2J_k^d}^2}{\mu(2J_k^d)} \int_{2J_k^d} V(y) \, d\mu(y) \geq C^{-1} a_k.}
The proof will be finished when we show that $a_k \to \8$. This follows from \eqref{asympt_2} by considering several cases. Let $a,b,c$ be such that $0\leq a<b<c$.

{\bf Case 1:} $J_k = [a,b]$, $I_k = [b,c]$.\\
{\bf Subcase 1:} $4a\leq 2b\leq c$. Using \eqref{asympt_2} we have
    \eqx{
    a_k \simeq (c/b)^2 (b/c)^{1+\a} = (c/b)^{1-\a} \simeq (|I_k|/|J_k|)^{1-\a}.
    }
{\bf Subcase 2:} $4a\leq 2b\geq c$. This subcase can hold only for finite $k$.\\
{\bf Subcase 3:} $4a\geq 2b\leq c$.
Using \eqref{asympt_2} we have
    \eqx{
    a_k \simeq \frac{c^2}{(b-a)^2} \frac{(b-a) {a^{\a}}}{c^{\a+1}} = \frac{{a^\a} c}{c^\a(b-a)} \geq \frac{c^{1-\a}}{(b-a)^{1-\a}} \simeq (|I_k|/|J_k|)^{1-\a}.
    }
{\bf Subcase 4:} $4a\geq 2b\geq c$.
Using \eqref{asympt_2} we have
    \eqx{
    a_k \simeq \eee{\frac{c-b}{b-a}}^2 \frac{(b-a) {a^{\a}}}{(c-b){b^{\a}}} \simeq |I_k|/|J_k|.
    }

{\bf Case 2:} $J_k = [b,c]$, $I_k = [a,b]$. Then
\eqx{
a_k \geq C \eee{\frac{|I_k|}{|J_k|}}^2 \frac{|J_k| b^\a}{|I_k| b^\a} \simeq |I_k|/|J_k|.
}
}

Recall that $\innprod{f,g} = \int_X f g \, d\mu$, so that $\sprod{-\Bb \phi, \psi}_\mu = \sprod{\phi, -\Bb \psi}_\mu$ for appropriate $\psi, \phi$. For $y>0$ the equation
    \eqx{
    \sprod{-\Bb \phi_y, \psi}_\mu = \psi(y)
    }
has the solution given by $\phi_y(x) = \frac{1}{2(1-\a)}\abs{x^{1-\a} - y^{1-\a}}$. We shall use $\phi_y$ to construct superharmonic functions that will be crucial in the proof of \eqref{eq:D}.

\lem{lem:conditionD}{
	Let $\a\in(0,1)$, $0\leq V\in L^1_{loc}(X,\mu)$, and $\ii(V)$ is as in Proposition \ref{construction}. Then
	\eqx{\lab{DDD}
    \int_X \Kt{2^n\abs{I}^2}(x,y) \, d\mu(x) \leq C 2^{-\frac{1-\a}{2} n}
    }
    for $y\in I^{**}$, $I\in \ii(V)$, and $n \geq 0$.
}

\pr{
Let $I$ be a dyadic interval such that
    \eqx{
    \frac{|2I|^2}{\mu(2I)}\int_{2I} V d\mu \leq 1, \qquad \frac{|2I^d|^2}{\mu(2I^d)}\int_{2I^d} V d\mu > 1.
    }
From \eqref{beta} we have $I^{**}\subset 2I$. By a continuity argument there exists $J$ such that $2I\subset J \subset 2I^d$ and $\frac{\abs{J}^2}{\mu(J)} \int_{J}V \, d\mu = 1$. Let $J=(a,b)$ and observe that $|J|\simeq |I|$. Define
	\eq{\label{def_phi}
      \phi_I(x) = 1 + \frac{1}{2(1-\a)} \int_J V(y)\abs{x^{1-\a}-y^{1-\a}} \, d\mu(y).
      }
Fix $z \in I^{**}$.  By \eqref{costam} and the doubling condition,
	\eq{\label{asdf}
    \phi_I(z)\leq 1 + C \sup_{y,y'\in J} |y^{1-\a} - (y')^{1-\a}| \int_J V d\mu \simeq C.
    }
Also, we claim that for $x\in X$,
    \eq{\label{asdf2}
    \phi_I(x) \simeq 1 +\frac{\mu(J)\abs{x^{1-\a}-z^{1-\a}}}{\abs{J}^{2}}.
    }
Indeed, if $x$ is such that $|x^{1-\a} - z^{1-\a}|\leq C |J|^2 \mu(J)^{-1}$, with $C$ large enough, this follows exactly as in \eqref{asdf}. In~the opposite case $|x^{1-\a} - z^{1-\a}|\geq C |J|^2 \mu(J)^{-1}$, we have $|x^{1-\a} - z^{1-\a}| \simeq |x^{1-\a} - y^{1-\a}|$ for $y\in J$ and the claim follows.

Now we proceed to a crucial argument that uses superharmonicity. Observe that formal calculation gives
\eqx{
 -\opL\phi_I(x) = -\opB\phi_I(x) - V(x)\phi_I(x) = V(x)(\mathbbm{1}_J(x) - \phi_I(x)) \leq 0
 }
and, consequently,
\eqx{\dt \Kt{t}\phi_I(z) = \int_X \Kt{t}(z,x)(-\opL\phi_I(x)) \, d\mu(x) \leq 0.}
This leads to
    \eq{\label{super}
    \Kt{t}\phi_I(z) \leq \phi_I(z), \qquad t>0.
    }
However, $\phi_I$ is not in $\Dom(\opL)$ (or even in $L^2(X,\mu)$), thus we provide a detailed proof of \eqref{super} in Appendix.

\newcommand{\rhoo}{\theta}
Denote
	\eqx{\rhoo(t) = \int_X \Kt{t}(z,x) \, d\mu(x), \quad t>0.}
Our goal is prove that, there exists $c_0$ such that for every $n\in\NN$
    \eq{\lab{induction}
    \rhoo(2^n|I|^2) \leq c_0 2^{-\frac{1-\a}{2} n}.
    }
This will follow by induction argument. Let $t>0$. By \eqref{FKK} and \eqref{eq:ptapprox1},
	\eq{\label{twot}
	\Kt{2t}(z,x) = \int_X \Kt{t}(z,y)\Kt{t}(y,x) \, d\mu(y) \leq C \rhoo(t) \mu(B(x,\sqrt{t}))^{-1}
	}
and
	\eqx{\label{est_m}
    \rhoo(2t) = \int_{\abs{x^{1-\a}-z^{1-\a}} < R} \Kt{2t}(z,x) \, d\mu(x) + \int_{\abs{x^{1-\a}-z^{1-\a}} > R} \Kt{2t}(z,x) \, d\mu(x) = A_1 + A_2.
    }
Here $R>0$ will be specified later. By \eqref{asdf2}, \eqref{super}, and \eqref{asdf} we have
	\spx{
     A_2 &\leq \frac{\abs{J}^2}{R\mu(J)} \int_{\abs{x^{1-\a}-z^{1-\a}} > R} \Kt{2t}(z,x) \frac{\mu(J)\abs{x^{1-\a}-z^{1-\a}}}{|J|^2} \, d\mu(x)\\
     &\leq C R^{-1}|J|^2 \mu(J)^{-1} \int  \Kt{2t}(z,x) \phi_I(x) \, d\mu(x)\\
     &\leq C R^{-1}|J|^2 \mu(J)^{-1}.
     }
To estimate $A_1$ we use \eqref{twot} and \eqref{ballmeasure},
    \spx{
    A_1 &\leq C \rhoo(t) \int_{\abs{x^{1-\a}-z^{1-\a}} < R} \mu(B(x,\sqrt{t}))^{-1} \, d\mu(x)\\
    &\leq C \rhoo(t) t^{-1/2} \int_{\abs{x^{1-\a}-z^{1-\a}} < R} \, dx.
    }

	{\bf Case A:} $z^{1-\a} \geq 2R$, then
	\eqx{\int_{\abs{x^{1-\a}-z^{1-\a}} < R} \, dx = (z^{1-\a} +R)^{\frac{1}{1-\a}} - (z^{1-\a} -R)^{\frac{1}{1-\a}}  \leq C R z^{\a}.}
	In this case
	\eq{\lab{case1}
    \rhoo(2t) \leq c_1 \left( \rhoo(t) t^{-1/2}  R z^{\a} + R^{-1}|J|^2\mu(J)^{-1}\right).
    }

    {\bf Case B:} $z^{1-\a} < 2R$,
	\eqx{\int_{\abs{x^{1-\a}-z^{1-\a}} < R} \, dx  \leq C (z^{1-\a} +R)^{\frac{1}{1-\a}} \leq C R^{\frac{1}{1-\a}}.}
	In this case
	\eq{\lab{case2}
     \rhoo(2t) \leq c_2 \left( \rhoo(t) t^{-1/2}  R^{\frac{1}{1-\a}}  + R^{-1}|J|^2\mu(J)^{-1}\right).
     }

Now we proceed to the proof of \eqref{induction}. The first step, $\rhoo(|I|^2)\leq C$, follows simply by \eqref{FKK}. Assume that \eqref{induction} holds for some $n$. The proof is finished by considering four cases. Since the calculations are similar in all the cases we present a detailed argument only in Subcase 1.1.

\noindent
{\bf Case 1:} $\rho(0,I) \geq 2|I|$. In this case $\mu(I) \simeq |I| z^\a$.\\
{\bf Subcase 1.1:} $\rhoo(2^n|I|^2) \geq 2^{n/2}|I|^2 z^{-2}$. Observe that
    \spx{
    R_1 & := 2^{-1} \rhoo(2^n |I|^2)^{-1/2} (2^n |I|^2)^{1/4}|I|^{1/2} z^{-\a} \\
    & \leq 2^{-1} 2^{-n/4}|I|^{-1} z \, 2^{n/4} |I|^{1/2}|I|^{1/2} z^{-\a} \\
    & = z^{1-\a}/2.
    }
Therefore we can use \eqref{case1} with $R=R_1$ and $t=2^n|I|^2$ together with the induction hypothesis, getting
    \spx{
    \rhoo(2^{n+1}|I|^2) & \leq c_1[\rhoo(2^n|I|^2)(2^n|I|^2)^{-1/2}2^{-1} \rhoo(2^n |I|^2)^{-1/2} (2^n |I|^2)^{1/4}|I|^{1/2} z^{-\a} z^\a \\
    & + 2 \rhoo(2^n |I|^2)^{1/2} (2^n |I|^2)^{-1/4}|I|^{-1/2} z^{\a}|I|^2\mu(I)^{-1}]\\
    & \leq c_3 (\rhoo(2^n |I|^2)^{1/2} 2^{-n/4} + \rhoo(2^n |I|^2)^{1/2} 2^{-n/4} |I|^{-1} z^{\a}|I|^2(z^\a|I|)^{-1})  \\
    & = 2c_3 \rhoo(2^n|I|^2)^{1/2}2^{-n/4} \\
    &\leq 2c_3c_0^{1/2} 2^{-n\frac{1-\a}{4}} 2^{-\frac{n}{4}} \leq c_0 2^{-(n+1)\frac{1-\a}{2}}.
    }
The last inequality holds if we choose $c_0$ such that $c_0\geq (2c_3)^2 2^{1-\a}$.\\


{\bf Subcase 1.2:} $\rhoo(2^n|I|^2) \leq 2^{n/2}|I|^2 z^{-2}$. One easily checks that
    \eqx{
    R_2 := 2^{-1}\left( \rhoo(2^n |I|^2)^{-1} (2^n |I|^2)^{1/2}|I| z^{-\a}\right)^{\frac{1-\a}{2-\a}} > z^{1-\a}/2.
    }
Putting $R=R_2$ and $t=2^n|I|^2$ into \eqref{case2} and using the induction hypothesis,
    \spx{
    \rhoo(2^{n+1}|I|^2) &\leq c_4 \left( \rhoo(2^n|I|^2) (2^n|I|^2)^{-1/2}|I|^{\frac{1}{1-\a}}z^{-\frac{\a}{1-\a}}\right)^{\frac{1-\a}{2-\a}} \\
    &\leq c_4 \left( \rhoo(2^n|I|^2) (2^n|I|^2)^{-1/2}|I|\right)^{\frac{1-\a}{2-\a}} \\
    &\leq c_4 \left(c_0(2^{-n})^{\frac{1-\a}{2} +\frac{1}{2}}\right)^{\frac{1-\a}{2-\a}}\leq c_0 2^{-(n+1)\frac{1-\a}{2}}.
    }
In the last inequality we choose $c_0$ such that $c_0\geq c_4^{2-\a} 2^{(1-\a)(2-\a)/2}$. Notice that we have used $z\geq |I|$, which follows from $\rho(0,I) \geq 2|I|$.

\noindent
{\bf Case 2:} $\rho(0, I) \leq 2\abs{I}$. In this case $\mu(I) \simeq |I|^{\a+1}$. Notice that $z\leq 4 |I|$. \\
{\bf Subcase 2.1} $\rhoo(2^n|I|^2) \geq 2^{n/2}|I|^{2-\a} z^{\a-2}$. Observe that
    \eqx{
    R_3 := 2^{-1} \rhoo(2^n |I|^2)^{-1/2} (2^n |I|^{2})^{1/4} |I|^{\frac{1-\a}{2}} z^{-\a/2} < z^{1-\a}/2.
    }
Putting $R=R_3$ and $t=2^n|I|^2$ into \eqref{case1} and using the induction hypothesis,
    \spx{
    \rhoo(2^{n+1}|I|^2) &\leq c_5 \rhoo(2^n|I|^2)^{1/2}2^{-n/4} z^{\a/2} |I|^{-\a/2} \\
    &\leq c_5 c_0^{1/2} 2^{-n\frac{1-\a}{4}} 2^{-\frac{n}{4}} z^{\a/2} |I|^{-\a/2}\\
    &\leq c_5 c_0^{1/2} 2^{\a} 2^{-n\frac{1-\a}{4}} 2^{-\frac{n}{4}}\\
    &\leq c_0 2^{-(n+1)\frac{1-\a}{2}}.
    }
The last inequality holds if we choose $c_0$ such that $c_0\geq c_5^2 2^{1+\a}.$\\
{\bf Subcase 2.2} $\rhoo(2^n|I|^2) < 2^{n/2}|I|^{2-\a} z^{\a-2}$. One easily checks that
\eqx{
    R_4 := 2^{-1} \left( \rhoo(2^n |I|^2)^{-1} (2^n |I|^2)^{1/2}|I|^{1-\a}\right)^{\frac{1-\a}{2-\a}} > z^{1-\a}/2.
    }
Putting $R_4$ and $t=2^n|I|^2$ into \eqref{case2} we obtain
    \spx{
    \rhoo(2^{n+1}|I|^2) &\leq c_6 \left( \rhoo(2^n|I|^2) (2^n|I|^2)^{-1/2}|I|\right)^{\frac{1-\a}{2-\a}} \\
    &\leq c_6 \left(c_0(2^{-n})^{\frac{1-\a}{2} +\frac{1}{2}}\right)^{\frac{1-\a}{2-\a}}\leq c_0 2^{-(n+1)\frac{1-\a}{2}},
    }
similarly as in subcase 1.2.
}

\lem{lem:conditionK}{
	Let $\a\in(0,1)$, $0\leq V\in L^1_{loc}(X,\mu)$, and $\ii(V)$ is as in Proposition \ref{construction}. Then the pair $(V, \ii(V))$ satisfies \eqref{eq:K}.
}

	\pr{
{\bf Case 1:} $\rho(0, I) \leq 2\abs{I}$. In this case $\mu(I) \simeq |I|^{1+\a}$. Using \eqref{eq:ptapprox1}, \eqref{ballmeasure} and \eqref{beta},
    	\alx{
    			\int_0^{2t}\int_{I^{***}} P_s(x,y) V(y)\, d\mu(y) \, ds \leq C \int_0^{2t} s^{-\frac{1+\a}{2}} ds \cdot \int_{I^{***}} V d\mu \leq C t^{\frac{1-\a}{2}} \frac{\mu(I)}{|I|^2}\leq C \eee{\frac{t}{|I|^2}}^{\frac{1-\a}{2}}.
		    }
   		{\bf Case 2:}  $\rho(0, I) > 2\abs{I}$. In this case $\mu(I) \simeq |I|c_I^{\a}$, where $c_I$ denotes the center of $I$. Using \eqref{eq:ptapprox1}, \eqref{beta} and the doubling condition,
   		\alx{
    			\int_0^{2t} \int_{I^{***}} P_s(x,y) V(y)\, d\mu(y) \, ds & \leq C \int_0^{2t}  \int_{I^{***}} \mu(B(y,\sqrt{s}))^{-1} V(y) d\mu(y) ds \\
    			& \leq C \int_0^{2t} s^{-\frac{1}{2}} c_I^{-\a} ds \cdot  \int_{I^{***}} V(y) d\mu(y) \\
    			& \leq C t^{\frac{1}{2}} \frac{\mu(I)}{|I|^2} c_I^{-\a} \leq C \eee{\frac{t}{|I|^2}}^{\frac{1}{2}}.
    		}
		}

Combining Lemmas \ref{lem:conditionD} and \ref{lem:conditionK} we obtain Theorem \ref{thm:maintheorem2}.


\section*{Appendix}\lab{App}

The goal of this Appendix is to give a precise proof of the formula \eqref{super}. In order to do this we need to pay careful attention to boundary conditions near zero, c.f. \eqref{DomQ} and \eqref{DomL}.  Recall that $J=(a,b)$, $0\leq a <b$, and $\int_J V d\mu = \mu(J)|J|^{-2} =: c_J>0$. Moreover, by \eqref{def_phi}, $\phi_I \in L^\8_{loc}[0,\8)$, $\abs{\phi_I(x) -\phi_I(0)} \leq C x^{1-\a}$ at $x=0$, and
\eq{\label{eq:eqq}
\phi_I'(x) =\frac{1}{2} x^{-\a} \begin{cases} -  c_J, & x\leq a\\
 \int_a^x V(y) \, d\mu(y) -  \int_x^b V(y) \, d\mu(y), & x\in(a,b)\\
  c_J , & x\geq b,
\end{cases}
}
\eq{\label{eq:eq}
-\phi_I''(x) - \frac{\a}{x} \phi_I'(x) = - \mathbbm{1}_J(x) V(x).
}
Since $\phi_I$ is not bounded in the infinity, we shall need some cut-off functions:
\eqx{
\eta_{n}\in C_c^\8(X), \quad \eta_{n} \rightarrow \mathbbm{1}_X, \quad \chi_{[0,n]}(x) \leq \eta_{n}(x) \leq \chi_{[0,n+1]} (x), \quad \norm{\eta_{n}'}_\8+\norm{\eta_{n}''}_\8 \leq C.
}
\lem{lem:lemat}{
	Let $n\in \NN$, $\psi\in \Dom(Q)$, $\psi\geq 0$, $\supp(\psi) \subseteq [0,n]$. Then $ \phi_I \eta_{n+1} \in \Dom(Q)$ and
\eq{\label{ineq}
	Q(\psi,\phi_I\eta_{n+1}) \geq 0.
}
}
\pr{
To simplify the notation we denote $\eta = \eta_{n+1}$. Let us first prove that $\phi_I\eta \in \rm{Dom}(Q)$. It is clear that $\sqrt{V} \phi_I\eta \in L^2(X, \mu)$, since $\rm{supp} \, \eta $
 is compact, $\phi_I, \eta \in{L^\8_{loc}(X)} $ and $V\in L^1_{loc}(X)$.  Observe that $\phi_I \eta \in C^1(X)$, but near zero we only have $\abs{(\phi_I \eta)'(x)} \leq C x^{-\a}$ . Let $\tau_k\in C_c^\8 (X)$ be such that
\eqx{
\chi_{[2/k,\8)}(x) \leq \tau_k(x) \leq \chi_{[1/k,\8)}(x) \ \text { and} \ \norm{\tau_k'}_\8 \leq Ck.
}
Define a sequence $\kappa_k =
\phi_I\eta \tau_k + \phi_I(0) \eta (1-\tau_k)$. It is clear that $\kappa_k \in  C_c^1[0,\8)$, $\kappa_k'(0^+) =0$, and $\kappa_k \to \phi_I \eta$ in $L^2(X,\mu)$. Moreover,
\eqx{
\kappa_k' = (\phi_I\eta)' \tau_k + \phi_I(0) \eta' (1- \tau_k) + ( \phi_I- \phi_I(0))\eta \tau_k' = (\phi_I\eta)' \tau_k + r_k.
}
Observe that $\eta' (1- \tau_k)=0$ and
\spx{
\norm{ ( \phi_I- \phi_I(0))\eta \tau_k'}_{L^2(X,\mu)}^2
\leq C \int_{1/k}^{2/k} x^{2(1-\a)} \norm{\tau_k'}^2_{\8} \, d\mu(x)\leq C k^{\a-1},
}
thus $\norm{r_k}_{L^2(X,\mu)} \to 0$. Now, $(\phi_I \eta - \kappa_n)' = (\phi_I\eta)' (1-\tau_k) -r_k \to 0$ in $L^2(X,\mu)$ since  $(\phi_I\eta)' \in L^2(X,\mu)$. This finishes the proof that $\phi_I\eta \in \rm{Dom}(Q)$.

What is left is to check \eqref{ineq}. Since $\phi_I \geq 1$, $\psi\geq 0$, and $-\phi_I'(x)x^\a|_{x=0^+} = c_J/2 >0$ the proof follows from a formal calculation based on \eqref{eq:eq}, namely
    \eq{\label{eq:equation}
    \int_X \psi'(x) (\phi_I(x) \eta(x))'\, d\mu(x)
    = -\int_X \psi(x) \mathbbm{1}_J(x) V(x)\, d\mu(x)  -\psi(x)\phi_I'(x)x^{\a}|_{x=0^+} .
    }
In the rest of the proof, for reader's convenience, we provide a detailed argument for \eqref{eq:equation}. The main problem is to deal with the boundary $x=0$.  Since $\psi \in \rm{Dom}(Q)$, we can find $\psi_k \in C_c^1[0,\8)$, $\psi'_k(0^+)=0$, such that $\psi_k \to \psi, \ \psi_k' \to \psi'$ in $L^2(X,\mu)$. We can additionally assume that $\supp(\psi_k)\subseteq [0,n+1]$.
By integrating by parts,
\alx{
\int_X \psi' (\phi_I \eta)' d\mu &= \lim_{k\to \8} \int_X \psi_k' \kappa_k' \, d\mu = \lim_{k\to \8} \int_X \psi_k' (\phi_I \eta)' \tau_k \, d\mu \\
&= - \lim_{k\to \8}\int_X \psi_k \eee{\eee{\phi_I \eta }'' + \frac{\a}{x}\eee{\phi_I \eta }'} \tau_k \, d\mu - \lim_{k\to \8} \int_X \psi_k \eee{\phi_I \eta}' \tau_k' \, d\mu\\
&= A_1 + A_2
}
Since, $\eta = 1$ on $\supp(\psi_k)$, by \eqref{eq:eq} we get $A_1 \to - \int_X \psi \mathbbm{1}_J V \, d\mu$ as $k\to \8$.
Notice that $ \eee{\phi_I \eta}'(x) = \phi_I'(x) = - x^{-\a} (c_J/2 - p(x))$ on $(1/k,2/k) $, where $p(x)=0$ when $a>0$ or $p(x)=\int_0^x V \,d\mu$ when $a=0$. We shall consider only the latter case, getting
\alx{A_2 &= \limk \int_{1/k}^{2/k} \psi_k \eee{\frac{c_J}{2} - p(x)} \tau_k' \, dx\\
&= \limk \eee{ \frac{ c_J}{2} \psi_k(2/k) - \psi_k(2/k) p(2/k)  -   \int_{1/k}^{2/k}  \psi_k' \eee{\frac{c_J}{2} - p(x)}  \tau_k \, dx + \int_{1/k}^{2/k} \psi_k V \tau_k \, d\mu(x)}\\
& = A_3+A_4+A_5+A_6. 
}
From the fact that $\psi,\psi_n, \psi', \psi_n' \in L^2(X)$ one can deduce that $\psi, \psi_n \in C[0,\8)$. Moreover, by the Cauchy-Schwarz inequality,
\eqx{
\abs{\psi_k(x)-\psi_k(0)} = \abs{\int_{0}^{x} \psi_k'(y) dy} \leq C \norm{\psi_k'}_{L^2(X,\mu)} x^{(1-\a)/2} \leq C x^{(1-\a)/2}.
}
As a consequence we get that $\psi_k(0) \to \psi(0)$ and $A_3 \to  \phi(0) c_J/2 > 0$. The proof is finished by noticing that $A_4,A_6 \to 0$ and
\spx{
|A_5| &\leq C \int_{1/k}^{2/k} \psi_k'(x) dx \leq C \sup_j \norm{\psi_j'}_{L^2(X,\mu)} k^{(\a-1)/2} \to 0.
}
}

Recall that $z\in J$ is fixed and denote
$
	\vt(u) = \Kt{u} \phi_I(z).
$
Our goal is to prove that $\vt(t+s)\leq \vt(t)$ for $t,s>0$.
Denote $k(x):=\Kt{t}(x,z) = \Kt{t/2}(\Kt{t/2}(\cdot,z))(x)$. Since the semigroup $\Kt{t}$ is analytic and $\Kt{t/2}(\cdot,z)\in L^2(X,\mu)$ (see \eqref{FKK} and \eqref{eq:ptapprox1}), we have $k\in \Dom(\opL)\subset \Dom(Q)$.


First, observe that
\alx{
	(\Kt{u}(k)\eta_{n-1})'(\phi_I\eta_{n+1})' & = (\Kt{u}(k)\eta_{n-1})'\phi_I' =  \Kt{u}(k)'\phi_I'\eta_{n-1}+ \Kt{u}(k)\phi_I'\eta_{n-1}'.
}
Using this,
\sp{\label{spsp}
	Q(\Kt{u}(k),\phi_I\eta_{n-1}) & = \int_X \Kt{u}(k)'(\phi_I\eta_{n-1})' \, d\mu + \int_X V \Kt{u}(k) \phi_I \eta_{n-1} \, d\mu \\
	& = \int_X {\Kt{u}(k)'\phi'_I\eta_{n-1}} \, d\mu + \int_X \Kt{u}(k)'\phi_I \eta'_{n-1} \, d\mu + \int_X V \Kt{u}(k) \phi_I \eta_{n-1} \, d\mu \\
	& = \int_X (\Kt{u}(k)\eta_{n-1})'(\phi_I\eta_{n+1})' \, d\mu - \int_X \Kt{u}(k) \phi'_I \eta'_{n-1}\, d\mu + \int_X \Kt{u}(k)'\phi_I \eta'_{n-1} \, d\mu \\
& \quad + \int_X V (\Kt{u}(k) \eta_{n-1}) (\phi_I \eta_{n+1}) \, d\mu \\
	& = Q(\Kt{u}(k)\eta_{n-1}, \phi_I \eta_{n+1})\\
 & - 2\int_X \Kt{u}(k) \phi'_I \eta'_{n-1}\, d\mu - \int_X \Kt{u}(k) \phi_I \eta''_{n-1} \, d\mu - \int \Kt{u}(k) \phi_I \eta'_{n-1} \frac{\a}{x} \, d\mu\\
	& = Q(\Kt{u}(k)\eta_{n-1}, \phi_I \eta_{n+1}) - 2B_1 - B_2-B_3.
}

\prop{decreasing}{ For fixed $I$ and $z\in J\supseteq 2I$ the function ${\vt(u)= \Kt{u} \phi_I(z)}$ is non-increasing.}
\pr{ Observe that $\left(\Kt{s}(k) - k\right) \phi_I$ is in $L^1(X, \mu)$. Using \eqref{spsp},
\spx{
	\vt(t+s) - \vt(t)
	& = \int_X \left(\Kt{s}(k) - k\right)(x) \phi_I(x) \, d\mu(x) \\
	& = \limn \int_X (\Kt{s}(k) - k)(x) \phi_I(x) \eta_{n-1}(x) \, d\mu(x) \\
	& = \limn \int_X \left( \int_0^s (-\opL) \Kt{u}(k)(x) \, du \right) \phi_I(x) \eta_{n-1}(x) \, d\mu(x) \\
	& = - \limn \int_0^s \int_X \opL \Kt{u}(k)(x) \phi_I(x) \eta_{n-1}(x) \, d\mu(x) \, du \\
	& = - \limn \int_0^s Q(\Kt{u}(k),\phi_I\eta_{n-1}) \, du \\
	& = - \limn \int_0^s \left( Q(\Kt{u}(k)\eta_{n-1}, \phi_I \eta_{n+1}) - 2 B_1 - B_2 -B_3\right) \, du.
}
Having in mind Lemma \ref{lem:lemat} it is enough to show that $\int_0^s B_i \, du \to 0$ as $n\to \8$ for $i=1,2,3$. This follows from Lebesgue's Dominated Convergence Theorem and the estimates we have already established. For example, for $B_2$ observe that $\abs{\phi_I(x) \eta_{n-1}''(x)} \leq C  |x|^{1-\a}\mathbbm{1}_{{[n-1,n]}}(x)$ for $n\geq N$ with $N$ large enough. Then a majorant is
    \alx{
    \int_0^s \int_X &\sup_{n\geq N} \abs{\Kt{u}(k)(x) \phi_I(x) \eta_{n-1}''(x)}\, d\mu(x) \, du \\
    &\leq  C \int_0^s \int_X \frac{|x|^{1-\a}}{\mu(B(x,\sqrt{t+u}))} \exp{\eee{-\frac{|x-z|^2}{c_2(t+u)}}} x^\a\, dx \, du\\
    &\leq  C \int_0^s \int_X \frac{|x|^{1-\a}}{(t+u)^{1/2}} \left(\frac{x}{x+\sqrt{t+u}}\right)^{\a} \exp{\eee{-\frac{|x-z|^2}{c_2(t+u)}}}\, dx \, du\\
    &\leq  C t^{-1/2} s \int_X |x|^{1-\a} \exp{\eee{-\frac{|x-z|^2}{c_2(t+s)}}} \, dx \\
    &\leq C(I,t,s).
    }
The integrals with $B_3$ and $B_1$ goes similarly. For the latter one we use  {$\abs{\phi_I'(x)} \leq C x^{-\a}$.}
}

\noindent
{\bf Acknowledgments:} The authors would like to thank Jacek Dziuba\'nski and the referees for their helpful comments.

\bibliographystyle{mn}        

\providecommand{\WileyBibTextsc}{}
\let\textsc\WileyBibTextsc
\providecommand{\othercit}{}
\providecommand{\jr}[1]{#1}
\providecommand{\etal}{~et~al.}



\end{document}